%%%%%%%%%%%%%%%%%%%%
%%% This the file for the revised version of the paper,
%%% Principal $G$-bundles over elliptic curves, by
%%% Robert Friedman,  John W. Morgan and Edward Witten
%%% Typeset in AMS-TeX, with the amsppt style
%%%%%%%%%%%%%%%%%%%%%

\input amstex

\define\scrO{\Cal O}
\define\Pee{{\Bbb P}}
\define\Zee{{\Bbb Z}}
\define\Cee{{\Bbb C}}
\define\Ar{{\Bbb R}}
\define\Pic{\operatorname{Pic}}
\define\Ker{\operatorname{Ker}}
\define\Sym{\operatorname{Sym}}
\define\Hom{\operatorname{Hom}}
\define\Aut{\operatorname{Aut}}
\define\Ext{\operatorname{Ext}}

\define\ad{\operatorname{ad}}

\define\endstatement{\endproclaim}
\define\theorem#1{\proclaim{Theorem #1}}
\define\lemma#1{\proclaim{Lemma #1}}
\define\proposition#1{\proclaim{Proposition #1}}
\define\corollary#1{\proclaim{Corollary #1}}
\define\claim#1{\proclaim{Claim #1}}

\define\section#1{\specialhead #1 \endspecialhead}

\documentstyle{amsppt}

\topmatter
\title Principal $G$-bundles over elliptic curves
\endtitle
\author {Robert Friedman,  John W. Morgan and Edward Witten}
\endauthor
\address Department of Mathematics, Columbia University, New
York,  NY 10027, USA\endaddress
\email rf\@math.columbia.edu, jm\@math.columbia.edu  \endemail
\address School of Natural Sciences, Institute for Advanced
Study, Princeton NJ 08540 USA\endaddress
\email  witten\@ias.edu  \endemail
\thanks The first author was partially supported by NSF grant
DMS-96-22681. The second author was partially supported by NSF
grant DMS-94-02988. The third author was partially supported
by NSF grant PHY-95-13835.
\endthanks

\endtopmatter

\nologo

\document

\section{1. Introduction.}

Let $E$ be an elliptic curve with origin $p_0$, and let $G$ be
a complex simple algebraic group. For simplicity, we shall
only consider the case where $G$ is simply connected, although
all of the methods discussed below can be extended to the case
of a general group $G$. The goal of this note is to announce
some results concerning the moduli of principal holomorphic
$G$-bundles over $E$. Detailed proofs, as well as a more
thorough discussion of the case where $E$ is allowed to be
singular or to vary in families and of the connections with del
Pezzo surfaces, elliptic
$K3$ surfaces, and Calabi-Yau manifolds which are elliptic or
$K3$ fibrations, will appear elsewhere.  

Grothendieck \cite{21} considered principal holomorphic
$G$-bundles over $\Pee ^1$, and showed that it was always
possible to reduce the structure group to a Cartan subgroup,
i\.e\. to a maximal (algebraic) torus in $G$. Atiyah \cite{1}
classified all holomorphic vector bundles over an elliptic
curve (in other words, the cases
$G=SL(n,
\Cee)$ or $G =PGL(n,
\Cee)$), without however considering the problem of trying to
construct a moduli space or find a universal bundle. In
\cite{16}, \cite{17}, and \cite{18}, this problem is studied
in the rank two case with a view toward constructing relative
moduli spaces in families. This approach has been generalized
to arbitrary rank in \cite{20}. A great deal of work has
been done on the moduli spaces and stacks of
$G$-bundles over a curve of genus at least two, partly
motivated by the study of conformal blocks and the Verlinde
formulas, by very many authors, e\.g\.
\cite{5}, \cite{15}. A basic method here is to
relate the moduli stack to an appropriate loop group. Related
constructions in the case of genus one have been carried out
by Baranovsky-Ginsburg \cite{4}, based on unpublished work of
Looijenga (see for example \cite{13}). They relate semistable
$G$-bundles to conjugacy classes in a corresponding affine
Kac-Moody group.  Recently Br\"uchert
\cite{9} has constructed a Steinberg-type cross-section for
the adjoint quotient of the affine Kac-Moody group whose image
lies in the set of regular elements, and this construction
leads to a moduli space for semistable
$G$-bundles which is equivalent to the one we construct in
Section 4 below. (We are indebted to Slodowy for calling our
attention to the work of Br\"uchert and sketching an argument
for the equivalence of the approach described above with the
one we give in this paper.) Finally, many of the results in
this note, along with applications to physics, are discussed in
\cite{19}.

The contents of this note are as follows. We will be concerned
with the classification of semistable $G$-bundles. As is
typical in invariant theory or moduli problems, the
classification will be up to a coarser equivalence than
isomorphism, which is usually called S-equivalence and will be
defined more precisely in Section 2. In Section 2, we describe
the moduli space of semistable
$G$-bundles over $E$ via flat connections for the maximal
compact subgroup $K$ of $G$, or equivalently via conjugacy
classes of representations $\rho\: \pi _1(E) \to K$. Such
bundles, which for a simply connected group
$G$ are exactly the bundles whose structure group reduces to a
Cartan subgroup, have an automorphism group which is as large
as possible in a certain sense within a fixed S-equivalence
class. The main result here is a theorem due to Looijenga and
Bernshtein-Shvartsman which describes this moduli space as a
weighted projective space. At the end of the section, we
connect this description, in the case where
$G=E_6, E_7, E_8$, with the moduli space of del Pezzo surfaces
of degree $3,2,1$ respectively and with the deformation theory
of simple elliptic singularities. In Section 3, we describe
regular $G$-bundles, which by contrast with flat bundles have
automorphism groups whose dimensions are as small as
possible within a fixed S-equivalence
class.  The generic $G$-bundle is both flat and regular.
However at special points of the moduli space we can choose
either a unique flat representative or a unique regular
representative, and it is the regular representatives which
fit together to give holomorphic families. In Section 4, we
show how special unstable bundles over certain maximal
parabolic subgroups can be used to give another description of
the moduli space in terms of regular bundles and obtain a new
proof of the theorem of Looijenga and Bernshtein-Shvartsman.
Finally, in the last section we discuss the existence of
universal bundles and give a brief description of how our
construction can be twisted with the help of a certain
spectral cover.

\section{2. Split semistable bundles.}

We fix notation for the rest of this paper. As before, $E$
denotes an elliptic curve with origin $p_0$. Let $G$ be a
simple and simply connected complex Lie group of rank $r$, and
let $\xi \to E$ be a holomorphic principal $G$-bundle over
$E$. The following definition differs from that given in
Ramanathan
\cite{32}, but is equivalent to it.

\definition{Definition 2.1} The principal bundle $\xi \to E$
is {\sl semistable\/} if the associated vector bundle $\ad
\xi$ is a semistable vector bundle. The principal bundle $\xi
\to E$ is {\sl unstable\/} if it is not semistable.
\enddefinition

Note that, if $\xi$ is stable in the sense of \cite{32}, it is
still possible for the vector bundle
$\ad
\xi$ to be strictly semistable. However, in our case ($G$
simply connected), there are essentially no properly stable
bundles over $E$, and so the above definition will suffice for
our purposes.

If $\xi$ is an unstable bundle, the structure group of $\xi$
reduces canonically to a parabolic subgroup $P$ of $G$, the
{\sl Harder-Narasimhan parabolic\/} associated to $\xi$ (see
for example \cite{31} or \cite{2}, pp\. 589--590). The
canonical reduction holds over a general base curve. In the
case of a base curve $E$ of genus one, it is easy to see that
the structure group further reduces to a Levi factor of $P$. 

Recall the following standard terminology: a {\sl family\/} of
principal 
$G$-bundles over $E$ parametrized by a complex space (or
scheme) $S$ is a principal $G$-bundle
$\Xi$ over $E\times S$. The family $\Xi$ is a family of {\sl
semistable\/}   principal $G$-bundles over $E$ if $\Xi|E\times
\{s\}=\Xi _s$ is semistable for all $s\in S$.  Finally, let
$\xi$ and
$\xi'$ be two semistable bundles over
$E$. We say that
$\xi$ and
$\xi'$ are {\sl S-equivalent\/} if there exists a family of
semistable bundles
$\Xi$ parametrized by  an irreducible $S$ and a point $s\in S$
such that, for
$t\neq s$, $\Xi |E \times \{t\} \cong \xi$ and $\Xi |E \times
\{s\} \cong
\xi'$. More generally, we let S-equivalence be the equivalence
relation generated by the above relation.

The following holds only under our assumption that $G$ is
simply connected.

\proposition{2.2} Let $\xi$ be a semistable principal
$G$-bundle, and suppose that the rank of $G$ is $r$. Then
$h^0(E; \ad \xi) \geq r$. Equivalently, $\dim \Aut _G\xi \geq
r$, where $\dim \Aut _G\xi$ denotes the group of global
automorphisms of $\xi$ \rom(as a
$G$-bundle\rom).
\endstatement

\definition{Definition 2.3} Let $\xi$ be a semistable
principal $G$-bundle. We call $\xi$ {\sl regular\/} if $h^0(E;
\ad \xi) = r$, or equivalently if $\dim
\Aut _G(\xi) = r$. We call $\xi$ {\sl split\/} if its
structure group reduces to a Cartan subgroup of
$G$, i\.e\. a maximal (algebraic) torus.
\enddefinition

It is easy to check that split bundles have the following
closure property: if there exists a family of semistable
bundles
$\Xi$ parametrized by  an irreducible $S$ and a point $s\in S$
such that, for
$t\neq s$, the bundles $\Xi |E \times \{t\}$ are split and all
isomorphic to each other, then
$\Xi |E \times \{s\}$ is isomorphic to  $\Xi |E \times \{t\}$,
$t\neq s$, and thus it is split as well. In general, however,
the condition of being split is neither open nor closed. On the
other hand, by the upper semicontinuity theorem, regularity is
an open condition: if
$\Xi$ is a family of semistable bundles parametrized by $S$
and $\Xi |E \times
\{s\}$ is regular, then $\Xi |E \times \{t\}$ is regular for
all $t$ in an open neighborhood of $s$.

To describe the set of split bundles, we introduce flat
bundles on the compact group. Let $K$ be a maximal compact
subgroup of $G$. Then $K$ is a compact, simple and simply
connected Lie group. If $\frak k$ is the Lie algebra of $K$ and
$\frak g$ is the Lie algebra of $G$, then $\frak g$ is the
complexification of
$\frak k$. Given a representation $\rho\: \pi _1(E) \cong \Zee
\oplus \Zee \to K$, we can form the associated principal
$K$-bundle $(\tilde E \times K)/\pi _1(E) \to E$, where $\pi
_1(E)$ acts on $\tilde E$, the universal cover of $E$, in the
usual way, and on $K$ via $\rho$. We shall call such a
$K$-bundle a {\sl flat\/} $K$-bundle. Using the inclusion
$K\subset G$, we can also view a flat $K$-bundle as a
$G$-bundle, and we shall also incorrectly refer to the induced
$G$-bundle as a flat $K$-bundle.  We will need the following
version of the theorem of Narasimhan-Seshadri
\cite{29} and Ramanathan
\cite{32} (see also Atiyah-Bott \cite{2} and Donaldson
\cite{12}):

\theorem{2.4} Let $\xi\to E$ be a semistable principal
$G$-bundle. Then there is a flat $K$-bundle $S$-equivalent to
$\xi$, and it is  unique  up to
isomorphism of flat $K$-bundles. More precisely, there is a
family of semistable principal $G$-bundles $\Xi$ over
$E\times \Cee$, such that, for
$t\neq 0$, $\xi_t=\Xi|E\times \{t\}\cong \xi$, and such that
$\xi _0= \Xi|E\times \{0\}$ is the $G$-bundle associated to a
flat $K$-bundle via the inclusion $K\subset G$. Finally, two
flat
$K$-bundles are isomorphic as
$G$-bundles if and only if they are isomorphic as $K$-bundles.
\endstatement

We note that Theorem 2.4 also holds for a non-simply connected
group. The special feature of simply connected groups which we
need to describe the moduli space of flat
$K$-bundles is contained in the following result of Borel
\cite{7} (see also
\cite{22} for the analogous algebraic result, due to Springer
and Steinberg):

\theorem{2.5} Let $K$ be a compact, simple, and simply
connected Lie group, and let $r_1$ and $r_2$ be two commuting
elements of $K$. Then there exists a maximal torus $T$ in $K$
with $r_1, r_2\in T$.
\endstatement

Since $\pi _1(E) \cong \Zee\oplus \Zee$, to give a
representation $\rho\: \pi _1(E) \to K$ is to give two
commuting elements $r_1, r_2\in K$. Thus a flat $K$-bundle
reduces to a $T$-bundle. In particular, we see that for a
simply connected group
$G$, every $G$-bundle associated to a flat bundle is split, and
conversely. On the other hand, if $G$ is not simply connected,
every split bundle lifts to the universal cover $\tilde G$ of
$G$, so that a $G$-bundle which does not lift to $\tilde G$
cannot be split.  Thus the correct notion for unliftable
bundles is that of a flat bundle.

Returning to the case of a simply connected group $G$, let $T$
be a maximal torus in the compact group $K$. One checks that
two homomorphisms from $\pi _1(E)$ to $T$ are conjugate by an
element of $K$ if and only if they are conjugate by an element
of the normalizer of $T$ in $K$. Thus we have:

\theorem{2.6} There is a natural bijection from the set of
flat $K$-bundles up to isomorphism, or equivalently the set of
semistable $G$-bundles up to S-equivalence, to the set $\Hom
(\pi _1(E), T)/W$, where $W$ is the Weyl group of $K$, acting
in the usual way on the maximal torus $T$.
\endstatement

Fix a maximal torus $T$ in
$K$. If $\Lambda = \pi _1(T)$, then $T\cong U(1) \otimes
_{\Zee}\Lambda$. Moreover, since $K$ is simply connected, if
$\frak t_\Ar$ denotes the real Lie algebra of $T$, then 
$\Lambda\subset
\frak t_\Ar$ is the lattice generated by the coroots $\alpha
\spcheck$, where
$\alpha \in \frak t_\Ar^*$ is a root. Now given a
homomorphism
$\rho\: \pi _1(E) \cong \Zee\oplus \Zee\to K$, the image of
$\rho$ is generated by two commuting elements of $K$ and so,
after conjugation, lies in $T$. The set of flat $T$-bundles is
naturally
$$\Hom (\pi _1(E), T) = \Hom (\pi _1(E),U(1) \otimes
_{\Zee}\Lambda ) \cong
\Hom (\pi _1(E), U(1))\otimes _\Zee \Lambda.$$ Now $\Hom (\pi
_1(E), U(1))$ is the set of flat line bundles on $E$, and is
naturally identified with $\Pic ^0E$. Since we have fixed a
base point $p_0\in E$, we can further identify $\Pic ^0E$ with
$E$. Thus the space of flat
$T$-bundles is naturally $E\otimes _\Zee \Lambda$. On the
other hand, as we are classifying not flat $T$-bundles but
flat $K$-bundles, we must take the quotient of $E\otimes _\Zee
\Lambda$ by the action of the Weyl group $W$ of $G$ acting on
$E\otimes _\Zee \Lambda$ via the natural action of $W$ on
$\Lambda$. We have thus described the coarse moduli space
of semistable $G$-bundles over $E$ as $(E\otimes _\Zee
\Lambda)/W$. A different proof of this result has been
given by Laszlo \cite{25}.

The varieties $(E\otimes _\Zee
\Lambda)/W$ have been studied by
Looijenga
\cite{27} and Bernshtein-Shvartsman
\cite{6}, who proved the following theorem:

\theorem{2.7} Let $E$ be an elliptic curve and let $\Lambda$
be the coroot lattice of a simple root system $R$ with Weyl
group $W$. Then $(E\otimes _\Zee
\Lambda)/W$ is a weighted projective space $WP(g_0, \dots,
g_r)$, where the weights $g_i$ are given as follows: $g_0 =
1$, and the remaining roots $g_i$ are found by choosing a set
of simple roots $\alpha _1, \dots , \alpha _r$, and then
writing the coroot $\tilde \alpha
\spcheck$ dual to the highest root
$\tilde \alpha$ as a linear combination $\sum _ig_i\alpha
_i\spcheck$ of the coroots dual to the simple roots. In case
$R$ is simply laced, we can identify the dual coroot  $\alpha
\spcheck \in R\spcheck$ to $\alpha$ with
$\alpha$, and consequently the
$g_i$ are the coefficients of $\tilde \alpha$ in terms of the
basis
$\alpha _1, \dots , \alpha _r$.
\endstatement

The proof of \cite{27} and \cite{6} makes use of formal
theta functions for a complexified affine Weyl group. We shall
outline a different proof of (2.7) below.

Since it will be important to motivate the construction of 
Section 4, let us give Looijenga's reason for studying the
space $(E\otimes _\Zee\Lambda)/W$. Let $(X, x_0)$ be the germ
of a simple elliptic singularity whose minimal resolution has
a single exceptional component which is a smooth elliptic
curve $E$ with self-intersection $-3, -2$, or $-1$. These are
exactly the simple elliptic singularities which can be
realized as hypersurface singularities in $(\Cee ^3, 0)$, and
we shall refer to them as being of type
$\tilde E_6, \tilde E_7, \tilde E_8$ respectively. These
singularities are weighted cones over
$E$ corresponding to a line bundle $L$ on $E$ of degree
$3$, $2$, or $1$, and thus have a
$\Cee ^*$-action. Moreover $\Cee ^*$ also acts on the tangent
space to the deformations of $(X, x_0)$. The zero weight
directions (in other words those directions fixed by the
$\Cee^*$-action) correspond to deforming
$(X, x_0)$ in an equisingular family by deforming $E$. The
remaining weights are negative, and deformations in the
neagative weight space correspond  to deforming $(X, x_0)$ to
a rational double point (RDP) singularity or smoothing it. The
local action of
$\Cee ^*$ on the negative weight deformations may be
globalized, and the quotient corresponding to the singularity
$\tilde E_r$ is a weighted projective space
$WP(g_0, \dots, g_r)$, where the weights $g_i$ are those
defined above for the root system
$E_r$. On the other hand, by the general theory of negative
weight deformations of singularities with $\Cee ^*$-actions,
and in particular by work of Pinkham \cite{30}, Looijenga
\cite{26}, and later M\'erindol \cite{28}, the points of
this weighted projective space parametrize triples $(\bar S, D,
\varphi)$, where $\bar S$ is a generalized del Pezzo surface
of degree $9-r$ (i\.e\., $\bar S$ has at worst rational double
point singularities and the inverse of the dualizing sheaf
$K_{\bar S}$ is ample on $\bar S$, with $K_{\bar S}^2 = 9-r$),
$D \in |-K_{\bar S}|$ is a smooth divisor, not passing through
the singularities of
$\bar S$,  and $\varphi$ is an isomorphism from $D$ to the
fixed elliptic curve
$E$ such that $\varphi ^*L = N_{D/\bar S}$. The moduli of such
triples 
$(\bar S, D, \varphi)$ can be described directly in terms of
the defining equations for $\bar S$ and can also be checked
directly to be a weighted projective space with the correct
weights. (Similar but slightly more involved arguments also 
handle the case of degree $4$ and $5$, in which case the
singularity is a codimension two complete intersection, in the
case  of degree $4$, and the corresponding root system is
$D_5$, or a Pfaffian singularity in case the degree is $5$,
and the root system is $A_4$.)

Now an elementary Torelli-type theorem shows that the
pair
$(\bar S, D)$ (ignoring the extra structure of $\varphi$) is
determined by the homomorphism
$\psi_0\: H^2_0(S; \Zee) \to D$, where $S$ is the minimal
resolution of $\bar S$ and $H^2_0(S;
\Zee)$ is the orthogonal complement of
$[K_S]$ in $H^2(S; \Zee)$, given as follows: represent a class
$\lambda\in H^2_0(S; \Zee)$ by a holomorphic line bundle $L$
on $S$ such that $\deg ( L|D) = 0$, and define $\psi_0
(\lambda)$ to be the element $L|D\in \Pic ^0D \cong D$. But
$H^2_0(S; \Zee)$ is  isomorphic to the root lattice for the
corresponding root system $E_r$, and this isomorphism is
well-defined up to the action of the Weyl group. The choice of
the isomorphism $\varphi$ enables one to extend the map
$\psi_0$ to a map $\psi\: H^2(S;\Zee)/\Zee[D] \to E$,
essentially because on the fixed curve $E$ we can choose a
$(9-r)^{\text{th}}$ root of the line bundle $L$, and
conversely the choice of such a root fixes an isomorphism from
$D$ to $E$ which lines up $L$ with $N_{D/\bar S}$. Now
$H^2(S;\Zee)/\Zee[D]$ is dual to the coroot lattice $\Lambda$
of the root system $E_r$, and 
$\psi$ defines an element of $E\otimes _\Zee\Lambda$,
well-defined modulo the action of $W$. In this way, we have
identified $WP(g_0, \dots, g_r)$ with $(E\otimes
_\Zee\Lambda)/W$. Let $\bar S$ be the result of contracting
all of the curves on $S$ not meeting $D$. Thus $\bar S$ has
certain rational double point (RDP) singularities. Under the
identification of the moduli space of pairs
$(S,D)$ with the set of $\psi\:
\Lambda \spcheck
\to E$, it is not difficult to show that the RDP singularities
on $S$ correspond to homomorphisms $\psi$ such that there is a
sub-root lattice 
$\Lambda ' \subseteq \Ker \psi$. In fact, the maximal such
lattice $\Lambda '$ describes the type of the RDP singularities
on $\bar S$. Here the main point is to show, by a Riemann-Roch
argument, that if
$\gamma \in \Ker \psi$ with $\gamma ^2 = -2$, then $\pm
\gamma$ is represented by an effective curve on $S$ disjoint
from $D$, and thus gives a singular point on the surface
obtained by contracting all such curves. In this way, there is
a link between subgroups of $E_r$, $r = 6,7,8$, and
singularities of the corresponding del Pezzo surfaces.

\section{3. Regular bundles.}

Recall that, for a simply connected group $G$, the bundle $\xi$
is {\sl regular\/} if
$h^0(E; \ad \xi)$ is equal to the rank of $G$. We begin by
giving a detailed description of the set of regular bundles in
case
$G$ is  one of the classical groups. At the end of the section
we shall outline the general structure of regular bundles. Let
us give a preliminary definition:

\definition{Definition 3.1} Let $I_n$ be the vector bundle of
rank $n$ and trivial determinant on $E$ defined inductively as
follows: $I_1=\scrO_E$, and
$I_n$ is the unique nonsplit extension of $I_{n-1}$ by
$\scrO_E$. More generally, if $\lambda$ is a line bundle on
$E$ of degree zero, we define
$I_n(\lambda) = I_n\otimes \lambda$.
\enddefinition

An easy argument shows that the algebra $\Hom (I_n, I_n)$ is
isomorphic to
$\Cee[t]/(t^n)$, and in particular it is a commutative
unipotent
$\Cee$-algebra  of dimension $n$. 

If $V$ is an arbitrary semistable vector bundle of degree zero
over $E$ and
$\lambda$ is a line bundle of degree zero over $E$, let
$V_\lambda \subseteq V$ be the sum of all of the subbundles of
$V$ which  are filtered by a sequence of subbundles whose
successive quotients are isomorphic to $\lambda$. An easy
argument shows that $V_\lambda$ itself is the maximal such
subbundle with this property and that $V = \bigoplus _\lambda
V_\lambda$. A straightforward induction classifies the
possible $V_\lambda$ as a direct sum
$\bigoplus _jI_{k_j}(\lambda)$. From this, it is easy to check:

\proposition{3.2} Let $V$ be a semistable vector bundle over
$E$ with trivial determinant, i\.e\. $V$ is a  principal
$SL(n)$-bundle over $E$. If
$V\cong \bigoplus _{i=1}^rI_{d_i}(\lambda _i)$, where the
$\lambda _i$ are line bundles on $E$ of degree zero, such that
$\lambda _1^{d_1} \otimes \cdots
\otimes \lambda _r^{d_r} =\scrO_E$ and $\sum _id_i = n$, then 
$V$ is  regular if and only if $\lambda _i\neq \lambda _j$ for
all $i\neq j$.
\endstatement

To deal with the case of the symplectic or orthogonal group,
the main point is to decide when a bundle $V$ carries a
nondegenerate alternating or symmetric form. The crucial case
is that of $I_n$. In this case, we have the following:

\proposition{3.3} There exists a nondegenerate alternating
pairing on $I_n$ if and only if $n$ is even.  There exists a
nondegenerate symmetric pairing on 
$I_n$ if and only if $n$ is odd. In both cases, every two such
nondegenerate pairings on $I_n$ are conjugate under the action
of $\Aut I_n$. 
\endstatement

With this said, we can describe the regular symplectic
bundles. It is simplest to describe them via the standard
representation:

\proposition{3.4} Let $V$ be a vector bundle of rank $2n$ over
$E$ with a nondegenerate alternating form, and suppose that
the dimension of the group of symplectic automorphisms of $V$
is $n$. Then there exist positive integers
$d_i$ and nonnegative integers $a_j$, $0\leq j\leq 3$, with
$\sum _id_i +
\sum _ja_j = n$, such that $V$ is isomorphic to 
$$\bigoplus _i\left(I_{d_i}(\lambda _i) \oplus I_{d_i}(\lambda
_i^{-1})\right)\oplus I_{2a_0} \oplus I_{2a_1}(\eta _1)
\oplus  I_{2a_2}(\eta _2)
\oplus  I_{2a_3}(\eta _3),$$ where the $\lambda _i$ are line
bundles of degree zero, not of order two, such that, for all
$i\neq j$, $\lambda _i\neq \lambda _j^{\pm1}$, and $\eta_1,
\eta _2, \eta _3$ are the three distinct line bundles of order
two on $E$. Conversely, suppose that $V$ is a vector bundle as
given above. Then $V$ has a nondegenerate alternating form,
all such forms have a group of symplectic automorphisms of
dimension exactly $n$,  and every two nondegenerate
alternating forms on
$V$ are equivalent  under the action of
$\Aut V$.
\endstatement

In particular, we see that a regular symplectic bundle is
always a regular bundle in the sense of $SL(2n)$-bundles.

For $SO(2n)$ and $SO(2n+1)$, the situation is a little more
complicated for two reasons. First, we shall only consider
those bundles which can be lifted to $Spin (2n)$ or $Spin
(2n+1)$, but shall not describe here the actual choice of a
lifting. Secondly, because of (3.3), it turns out that a
regular $SO(n)$-bundle does not always give a regular
$SL(n)$-bundle.  

\proposition{3.5} Let $V$ be a vector bundle of rank $2n$ over
$E$ with a nondegenerate symmetric form, and suppose that the
dimension of the group of orthogonal automorphisms of $V$ is
$n$. Finally suppose that $V$ can be lifted to a principal
$Spin (2n)$-bundle. Then
$V$ is isomorphic to 
$$\bigoplus _i\left(I_{d_i}(\lambda _i) \oplus I_{d_i}(\lambda
_i^{-1})\right)\oplus \bigoplus _j\left(I_{2a_j+1}(\eta
_j)\oplus\eta _j\right)$$ where the $\lambda _i$ are line
bundles of degree zero, not of order two, such that, for all
$i\neq j$, $\lambda _i\neq \lambda _j^{\pm1}$,  $\eta _0 =
\scrO_E, \eta_1,
\eta _2, \eta _3$ are the four distinct line bundles of order
two on $E$, and the second sum is over some subset
\rom(possibly empty\rom) of $\{0,1,2,3\}$. Conversely, every
such vector bundle
$V$ has a nondegenerate symmetric form, all such forms have a
group of orthogonal automorphisms of dimension exactly $n$,
and every two nondegenerate symmetric forms on $V$ are
equivalent  under the action of $\Aut V$.
\endstatement

Here the symmetric form on $I_{2a_0+1}\oplus\scrO_E$ consists
of the orthogonal direct sum of the  nondegenerate form on the
factor $I_{2a_0+1}$ given by (3.3), together with the obvious
form on $\scrO_E$, and similarly for the summands
$I_{2a_i+1}(\eta _i)\oplus \eta _i$. Moreover, not all of the
summands $I_{2a_j+1}(\eta _j)\oplus\eta _j$ need be present in
$V$.  We remark that, if a vector bundle
$\bigoplus _jI_{d_j}(\lambda _j)$ is isomorphic to its dual,
and the sum of all the factors where $\lambda _j =
\eta _i$ for some $i$ has odd rank, then the same must be true
for all of the
$\eta _i$. Thus, if the automorphism group of $V$ is to be as
small as possible, then either $V$ is as described in (3.5) or
$V$ is of the form 
$$\bigoplus _i\left(I_{d_i}(\lambda _i) \oplus I_{d_i}(\lambda
_i^{-1})\right)\oplus I_{2a_0+1}\oplus I_{2a_1+1}(\eta
_1)\oplus  I_{2a_2+1}(\eta_2)\oplus  I_{2a_3+1}(\eta _3).$$
But in this last case $V$ does not lift to a $Spin
(2n)$-bundle.

The case of $SO(2n+1)$, which we shall not state explicitly,
is completely analogous, except that the summand
$I_{2a_0+1}\oplus\scrO_E$ is replaced by the odd rank summand 
$I_{2a_0+1}$, which must always be present.

We return now to the study of regular bundles over a general
group $G$.

\proposition{3.6} Let $\xi$ be a semistable principal
$G$-bundle over $E$. Then the structure group of $\xi$ reduces
to an abelian subgroup of $G$. If furthermore $\xi$
is regular, the structure group of $\xi$ reduces to an abelian
subgroup of
$\Aut _G\xi$, which naturally sits inside $G$
up to conjugation.
\endstatement

In fact, one can take the structure group of $\xi$ to be of
the following form. Let $\xi _0$ be the split bundle
S-equivalent to $\xi$, corresponding to the
representation $\rho\: \pi _1(E) \to T\subset K$. Let $T_0$
be the image of $\rho$. Then there exists a
subgroup $U$ of $G$ commuting with $T_0$, which
is either trivial or a $1$-parameter commutative unipotent
subgroup, such that the structure group of $\xi$ reduces to
$T_0U$. 

We now describe
the set of bundles which are simultaneously regular and
split.  If $\xi$ is split, then
$\xi$ corresponds to a point of $(E\otimes_{\Zee} \Lambda)/W$.
After lifting this point to an element $\mu$ of
$E\otimes_{\Zee} \Lambda$, we see that we can describe
$\ad \xi$ as follows. A root $\alpha$ defines a homomorphism
$\Lambda \to
\Zee$, and thus a homomorphism $E\otimes _{\Zee}\Lambda \to E
\cong \Pic ^0E$. Denote the image of $\mu$ in $E$ by
$\alpha(\mu)$ and the corresponding line bundle by
$\lambda_{\alpha(\mu)}$. Then, as vector bundles,
$$\ad \xi \cong \scrO_E^r \oplus \bigoplus
_\alpha\lambda_{\alpha(\mu)}.$$ Hence $\xi$ is regular if and
only if, for every root $\alpha$, $\alpha (\mu)
\neq 0$.

In particular, there is a nonempty Zariski open subset of
$(E\otimes_{\Zee}
\Lambda)/W$ such that all of the corresponding split bundles
are regular. In fact, on this open subset, S-equivalence is
the same as isomorphism.

At the other extreme, we can consider bundles which are
S-equivalent to the trivial bundle. The split representative
for the S-equivalence class corresponds to the image of $0\in
E\otimes_{\Zee} \Lambda$ in $(E\otimes_{\Zee}
\Lambda)/W$, which has the unique preimage $0\in
E\otimes_{\Zee} \Lambda$. To describe the regular
representative, or more precisely its adjoint bundle, we first
recall  the definition of the {\sl Casimir weights\/}
$d_1, \dots, d_r$ of a root system
$R$. These can be defined to be the numbers $m_i +1$, where
the $m_i$ are the exponents of $R$ (cf\.
\cite{8}, V (6.2)), and they are also the degrees of a set of
homogeneous generators for the invariants of the symmetric
algebra of the vector space corresponding to the root system
$R$ under the action of the Weyl group. To describe a  regular
bundle S-equivalent to the trivial bundle, we shall describe
its adjoint bundle. (Here an $(\ad G)$-bundle has in general
finitely many liftings to a $G$-bundle, but exactly one of
these will turn out to be S-equivalent to the trivial bundle.)

\proposition{3.7} There is a unique regular $G$-bundle $\xi$
S-equivalent to the trivial bundle. As vector bundles over $E$,
$$\ad \xi \cong \bigoplus _iI_{2d_i-1},$$
where the $d_i$ are the Casimir weights of the root system of
$G$.
\endstatement

The bundle $\ad \xi$ can be seen to be an ($\ad G$)-bundle as
follows: start with the bundle $I_3 = \Sym ^2I_2$. It is an
$SL(2)$ bundle which descends to an
$SO(3)$-bundle. Now there is a ``maximal" embedding of $SO(3)$
in $G$, unique up to conjugation. Thus there is a
representation $\rho$ of $SO(3)$ on the Lie algebra $\frak g$.
Under this representation $\frak g$ decomposes as a direct sum 
$$\frak g = \bigoplus _i\Sym ^{2d_i-2}(\Cee ^2),$$ where we
view $\Cee ^2$ as the standard representation of $SL(2)$, and
thus its odd symmetric powers give representations of $SO(3)$.
In particular, the
$G$-bundle induced by $\rho$ gives rise to the ($\ad
G$)-bundle described above.

We can generalize the above picture for the trivial bundle to
an arbitrary bundle. Let $\xi$ be an arbitrary semistable
$G$-bundle and let $\xi _0$ be the unique split bundle
S-equivalent to $\xi$. Then $\Aut \xi _0$ is up to isogeny a
product of $N$  factors $G_i$, where each factor $G_i$ is
either simple or isomorphic to $\Cee^*$.  Let $\mu \in
E\otimes_{\Zee} \Lambda$ be a representative for the class of
$\xi _0$. The Lie algebra
$H^0(E;\ad \xi _0)$ of $\Aut \xi _0$ is identified with 
$$\frak h \oplus \bigoplus _{ \alpha(\mu) = 0}\frak g^\alpha,$$
where $\frak g^\alpha$ is the root space corresponding to the
root $\alpha$ (and thus in particular the rank of this
reductive Lie algebra is $r$). We then have:

\proposition{3.8} With notation as above, let
$\xi_{\text{reg}}$ be a regular semistable bundle S-equivalent
to $\xi_0$. Let $r_i$ be the rank of $G_i$, where by
definition $r_i =1$ if $G_i\cong \Cee^*$, and let $d_{ij},
1\leq j\leq r_i$ be the Casimir weights of $G_i$, where we set
$d_{i1} = 1$ if $G_i\cong
\Cee ^*$. Then the maximal subbundle of $\ad \xi_{\text{reg}}$
which is filtered by subbundles whose successive quotients are
$\scrO_E$ is
$$\left(\ad \xi_{\text{reg}}\right)_{\scrO_E} = \bigoplus
_{i=0}^N\bigoplus _{j=1}^{r_i}I_{2d_{ij}-1}.$$
\endstatement

From this, it is possible in principle to give a complete
description of $\ad
\xi_{\text{reg}}$.

As a consequence of Proposition 3.6, one can show:

\proposition{3.9} Let $\xi$ be a semistable principal
$G$-bundle. Then $\xi$ is S-equivalent to a unique regular
semistable bundle and to a unique split bundle. 
\endstatement

There are thus two canonical representatives for every
S-equivalence class, depending on whether we choose the
regular or the split representative. For an open dense subset
of bundles, these two representatives will in fact coincide.
As should be clear from Section 2, the split representatives
arise most naturally from the point of view of flat
connections. However, if we try to find a universal
holomorphic $G$-bundle, then we must work instead with regular
bundles. In fact, even working locally, it is not possible to
fit the split bundles together into a universal bundle, even
for $SL(n)$.

Finally, we make some comments about the automorphism group of
a regular bundle. 

\proposition{3.10} Let $\xi$ be a regular semistable
$G$-bundle. If $\Aut _G(\xi)$ is the automorphism group of
$\xi$ and $(\Aut_G(\xi))^0$ is the component of $\Aut
_G(\xi)$ containing the identity, then $(\Aut _G(\xi))^0$
is abelian. Moreover, $\Aut _G(\xi)$ is itself abelian if and
only if $\xi$ corresponds to a smooth point of the moduli
space of S-equivalence classes of semistable $G$-bundles.
\endstatement

In fact, a careful analysis of the root systems involved shows
that the singular locus of the moduli space corresponding to
$\Zee/d\Zee$-isotropy is smooth and irreducible, of dimension
equal to the number of $i$ such that $d|g_i$, in the notation
of Theorem 2.7. Of  course, this statement also follows
directly from Theorem 2.7.

\section{4. The parabolic construction.}

In this section, we describe a method of constructing families
of regular semi\-stable $G$-bundles. The motivation is as
follows: we seek to find an analogue for bundles of the
singularities picture outlined above in Section 2. That is, we
seek to find a mildly ``singular" (in other words, unstable)
$G$-bundle
$\xi _0$ together with a $\Cee ^*$-action on its deformation
space, such that the weighted projective space corresponding
to the quotient of the negative weight deformations of $\xi
_0$ by $\Cee^*$ is both the weighted projective space
$WP(g_0,
\dots, g_r)$ and is the coarse moduli space of semistable
$G$-bundles modulo S-equivalence. (Actually, with our
conventions the action of $\Cee^*$ will be by positive
weights.) It will also turn out that the points of the weighted
projective space parametrize regular
$G$-bundles, as opposed to split bundles, and will thus enable
us to find locally a universal
$G$-bundle away from the orbits where $\Cee ^*$ does not act
freely. In fact, in many cases we can use this construction to
produce a global universal $G$-bundle.

To pursue this idea further, we have seen that unstable
$G$-bundles over $E$ reduce to a parabolic subgroup of $E$,
and further to a Levi factor $L$. Conversely, fix a maximal
parabolic subgroup $P$ of $G$ and a Levi decomposition $P =
LU$, where $U$ is the unipotent radical of $P$ and $L$ is the
reductive or Levi factor. Then $U$ is normal, all Levi factors
are conjugate in $P$, and the quotient homomorphism
$P \to L$ is well-defined. The group $L$ is never semisimple;
in fact, since
$P$ is a maximal parabolic, the connected component of the
center of
$L$ is $\Cee ^*$. The maximal parabolic subgroup $P$ has a
canonical character $\chi\: P \to \Cee ^*$ (the unique
primitive dominant character), which is induced from a
character $L\to \Cee$. Using this character, we can define the
determinant line bundle of a principal
$L$-bundle over $E$. Fix  an
$L$-bundle
$\eta$, such that $\det \eta$ has negative degree. The induced
$G$-bundle $\xi _0$ is unstable, because $\xi_0$ also reduces
to the opposite parabolic to $P$, and the determinant line
bundle for the primitive dominant character of the opposite
parabolic has positive degree.  Consider the set of all
$P$-bundles
$\xi$ such that the associated
$L$-bundle (via the homomorphism $P\to L$) is $\eta$. It is
straightforward to classify all such bundles: the group $L$
acts by conjugation on $U$, and the
$L$-bundle
$\eta$ and the action of $L$ on $U$ define a sheaf of
unipotent groups
$\underline{U}(\eta)$ on $E$, which is in general nonabelian.
The set of all isomorphism classes of $P$-bundles
$\xi$ which reduce to $\eta$ may then be identified with the
cohomology set $H^1(E;
\underline{U}(\eta))$. The $\Cee^*$ in the center of $L$ then
acts on  $H^1(E; \underline{U}(\eta))$. Cohomology
sets similar to
$H^1(E; \underline{U}(\eta))$, arising from the $H^1$ of a
sheaf of unipotent groups over a base curve, have been studied
in a different context by Babbitt and Varadarajan
\cite{3}, following ideas of Deligne, as well as by Faltings
\cite{14}.  Using similar ideas, one can show that $H^1(E;
\underline{U}(\eta))$ has a (non-canonical) linear structure
and that $\Cee ^*$ acts linearly in this structure with
positive weights (following certain standard conventions), so
that the quotient is isomorphic to a  weighted projective
space. 

In the case of $SL(n)$, it is easy to make these ideas
explicit. The maximal parabolic subgroups of $SL(n)$ correspond
to filtrations $\{0\} \subset \Cee ^d \subset
\Cee ^n$, where $0< d < n$. For each such $d$, there is a
unique stable bundle
$W_d$ over $E$ of rank $d$ such that $\det W_d =
\scrO_E(p_0)$. The unstable bundle which we consider is then
$W_d^* \oplus W_{n-d}$, and it has a nontrivial
$\Cee^*$-action which acts trivially on $\det (W_d^* \oplus
W_{n-d})$. In this case, a straightforward argument shows:

\theorem{4.1} Let $V$ be a regular semistable vector bundle of
rank $n$. Then there is an exact sequence
$$0 \to W_d^* \to V \to W_{n-d} \to 0. $$ Moreover, the
automorphism group of $V$ acts transitively on the set of
subbundles of $V$ isomorphic to
$W_d^*$ whose quotients are isomorphic to $W_{n-d}$. Finally,
if $V$ is a nonsplit extension of $W_{n-d}$ by $W_d^*$, then
$V$ is in fact a regular  semistable vector bundle.
\endstatement

We note that in this case the parabolic subgroup in question is
$$P = \left\{\,\pmatrix A&B\\O&D\endpmatrix: A\in GL(d), D\in
GL(n-d), \det A\cdot
\det D = 1\,\right\},$$ the Levi factor of $P$ is given by
$$L= \left\{\,\pmatrix A&O\\O&D\endpmatrix: A\in GL(d), D\in
GL(n-d), \det A\cdot
\det D = 1\,\right\},$$ and the unipotent radical $U$ of $P$,
which in this case is abelian, is given by
$$U= \left\{\,\pmatrix I&B\\O&I\endpmatrix: B \text{ is a
$d\times (n-d)$ matrix}\,\right\}.$$ It is easy then to
identify $H^1(E; \underline{U}(\eta))$ with the usual sheaf
cohomology group $H^1(E; W_{n-d}^*\otimes W_d^*)$ and the
$\Cee ^*$-action with the usual one, up to a factor. In this
way, the moduli space of regular semistable vector bundles
over $E$ of rank $n$ and trivial determinant is identified
with $\Pee ^{n-1}$, a fact which could also be established by
spectral cover methods \cite{20}. The full
tangent space to the deformations of the unstable bundle
$W_d^* \oplus W_{n-d}$ keeping the determinant trivial is
$H^1(E; \ad(W_d^* \oplus W_{n-d}))$. This group contains the
subgroup $H^1(E; W_{n-d}^*\otimes W_d^*)$ which is tangent to
the set of extensions described above. The one remaining
direction has $\Cee^*$-weight zero, which correponds to moving
the point $p_0$ on $E$ and which should be viewed as a one
parameter family of locally trivial deformations.

In the case of $SL(n)$, or equivalently the root system
$A_{n-1}$, every maximal parabolic subgroup has an abelian
unipotent radical and there is an appropriate construction
from any such subgroup giving the moduli space of regular
semistable
$G$-bundles. In all other cases, we have the following:

\theorem{4.2} Let $G$ be a complex, simple, and simply
connected group, not of type $A_n$. Then there exists a unique
maximal parabolic subgroup $P$ of $G$, up to conjugation, such
that, if $L$ is the Levi factor of $P$, then there exists an
$L$-bundle $\eta$  with the following properties:
\roster
\item"{(i)}" The connected component of the automorphism group
of $\eta$ as an
$L$-bundle is $\Cee ^*$.
\item"{(ii)}" The line bundle $\det \eta$ has negative
degree, and so the
$G$-bundle
$\xi _0$ induced by
$\eta$ is unstable.
\item"{(iii)}" If $U$ is the unipotent radical of $P$, then
the nonabelian cohomology set $H^1(E; \underline{U}(\eta))$
has the structure of affine $(r+1)$-dimensional space.
\item"{(iv)}" There exists a linear structure on $H^1(E;
\underline{U}(\eta))$  for which the natural copy of $\Cee
^*\subseteq \Aut_G\xi _0$ acts linearly, fixing the trivial
element, and with negative weights. The
stabilizer of every nontrivial element of $H^1(E;
\underline{U}(\eta))$ is  finite, and the quotient $(H^1(E;
\underline{U}(\eta)) -
\{0\})/\Cee^*$ is a weighted projective space $WP(g_0, \dots,
g_r)$.
\item"{(v)}" If $\xi$ is a $P$-bundle over $E$ corresponding
to an element of 
$H^1(E; \underline{U}(\eta)) - 0$, then $\xi$ is a regular
semistable bundle.
\endroster In all cases, the bundle $\eta$ with the above
properties is uniquely specified by requiring that $\det \eta
= \scrO_E(-p_0)$. 
\endstatement

In fact, (iv) and (v) are a consequence of the other
properties. If we do not specify that $\det \eta =
\scrO_E(-p_0)$, then it is still the case that $\det\eta$ must
have degree $-1$, and so
$\eta$ is specified up to translation on $E$.

We note that all of the weights are equal, in other words the
weighted projective space is an ordinary projective space,
exactly in the cases $A_n$ and $C_n$, in other words for the
groups $SL(n+1)$ and $Sp(2n)$. In all other cases, for a simply
connected group $G$, the weighted projective space will in fact
have singularities.

To describe the maximal parabolic subgroups which arise in
Theorem 4.2, note first that maximal parabolic subgroups of
$G$, up to conjugation, are in one-to-one correspondence with
the vertices of the Dynkin diagram of the corresponding root
system. In case $G$ is 
$D_n$ or
$E_6, E_7, E_8$, the maximal parabolic subgroup in Theorem 4.2
corresponds to the unique trivalent vertex of the Dynkin
diagram. In the remaining cases, the vertex in question is the
unique vertex meeting the multiple edge which is the long
root. (Such vertices will be trivalent in an appropriate sense
except for the case $C_n$.)

Let us describe the construction explicitly for the remaining
classical groups. The simplest case after $A_n$ is the case of
$Sp(2n)$, in other words $C_n$. In this case the parabolic in
question corresponds to those elements of $Sp(2n)$ which
preserve a totally isotropic subspace of dimension $n$. Thus
$$P = \left\{\,\pmatrix T&B\\O&{}^tT^{-1}\endpmatrix: T\in
GL(n), T^{-1}B = {}^t(T^{-1}B)\,\right\},$$ the Levi factor of
$P$ is given by
$$L= \left\{\,\pmatrix T&O\\O&{}^tT^{-1}\endpmatrix: T\in
GL(n)
\,\right\},$$ and the unipotent radical $U$ of $P$, which in
this case is also abelian, is given by
$$U= \left\{\,\pmatrix I&B\\O&I\endpmatrix: {}^tB =
B\,\right\}.$$ The unstable symplectic bundle corresponding to
$\eta$ is the bundle
$W_n^*\oplus W_n$, with the first factor embedded as a totally
isotropic subbundle and the second as its dual. It is easy
then to identify
$H^1(E; \underline{U}(\eta))$ with the usual sheaf cohomology
$H^1(E; \Sym ^2W_n^*)$. Here $\Cee^*$ acts with constant
weight, so that the quotient is an ordinary (smooth)
$\Pee^{n-1}$.

Next we consider $Spin(2n)$, although here it will be more
convenient to work in $SO(2n)$. The natural analogue of the
construction for the symplectic group would lead to the
unstable bundle $W_n^*\oplus W_n$, together with the symmetric
nondegenerate form for which $W_n^*$ is isotropic and which
identifies the dual of $W_n^*$ with the complementary $W_n$.
Such orthogonal bundles do not lift to $Spin(2n)$, although
this construction does identify all of the regular semistable
$SO(2n)$-bundles with $w_2\neq 0$ with the projective space on
$H^1(E; \bigwedge ^2W_n^*)$, which is a $\Pee ^{n-2}$. For
liftable $SO(2n)$-bundles, we use the parabolic subgroup
corresponding to the trivalent vertex, which is the subgroup of
$g\in SO(2n)$ preserving an isotropic subspace of rank $n-2$.
In this case the unipotent radical is nonabelian. The bundle
$\eta$, viewed as an unstable
$SO(2n)$-bundle $\xi _0$, is the bundle
$$\xi _0 = W_{n-2}^* \oplus Q_4 \oplus W_{n-2},$$ where
$W_{n-2}^*$ is an isotropic subspace, $Q_4$ is the
$SO(4)$-bundle
$\scrO_E\oplus \eta _1\oplus \eta _2 \oplus \eta _3$, in the
notation of Section 3, with a diagonal nondegenerate symmetric
pairing, and
$Q_4$ is orthogonal to the direct sum $W_{n-2}^* \oplus
W_{n-2}$. More invariantly $Q_4 = Hom (W_2, W_2)$ with the
quadratic form given by the trace. Note that neither $Q_4$
nor  $W_{n-2}^* \oplus W_{n-2}$ lifts to a
$Spin$-bundle, and hence the direct sum is liftable. In this
case, the nonabelian cohomology set
$H^1(E; \underline{U}(\eta))$ (for $Spin(2n)$) has a weight
$1$ piece given by
$H^1(Q_4\otimes W_{n-2}^*)$, of rank $4$, and a weight $2$
piece given by
$H^1(E; \bigwedge ^2W_{n-2}^*)$, which has rank $n-3$. Similar
results hold for
$Spin (2n+1)$, by replacing $Q_4$ by 
$$Q_3 = \eta _1\oplus \eta _2 \oplus \eta _3 = \ad W_2.$$

Returning to the general case, let us show that the weighted
projective space
$WP(g_0, \dots, g_r)$ arising from the parabolic construction
can be naturally identified with $(E\otimes _\Zee\Lambda)/W$,
thus giving a new proof of Looijenga's theorem. One first
shows that there exists a universal $G$-bundle over
$E\times H^1(E; \underline{U}(\eta))$ in the appropriate
sense. By general properties, there is a $\Cee^*$-equivariant
map from $H^1(E; \underline{U}(\eta))-0$ to the moduli space
of semistable $G$-bundles, in other words to $(E\otimes
_\Zee\Lambda)/W$.

\theorem{4.3} The induced map $WP(g_0, \dots, g_r) \to
(E\otimes _\Zee\Lambda)/W$ is an isomorphism.
\endstatement 

The essential point of the proof is to compare
the determinant line bundles on the two sides, and then to use
the elementary fact that a degree one morphism from a weighted
projective space to a normal variety is an isomorphism. On the
weighted projective side, the determinant line bundle is
always Cartier, and in fact it is the line bundle
$K_{WP^r}^{2}$. On the other side, it is easy to calculate
the preimage of the determinant line bundle in $E\otimes
_\Zee\Lambda$. At least in the case of a simply laced root
system $R$, the fact that the degree of the morphism in
question is one then follows from the fact that the order of
the Weyl group is $r!(g_1\cdots g_r)\det R$ \cite{8}.

The parabolic construction also leads to a proof of the
existence of universal bundles in certain cases. For a fixed 
$G$, we denote by $\Cal M_E=\Cal M_E(G)$ the moduli space of
regular semistable $G$-bundles over $E$ and by $\Cal M_E^0$
the smooth locus of $\Cal M_E$.

\theorem{4.4} If $G= SL(n)$, let $P_d$ be the maximal
parabolic subgroup of $SL(n)$ stabilizing the flag $\{0\}
\subset \Cee^d\subset \Cee^n$, and if $G\neq SL(n)$, let
$P$ be the maximal parabolic subgroup of $G$ described in
Theorem \rom{4.2}. Let
$n_P$ be the positive  integer defined as follows\rom:
\roster
\item"{(i)}" If $G=SL(n)$ and $P=P_d$, then $n_{P_d} =
n/\gcd(d,n)$.
\item"{(ii)}" If $G$ is of type $C_n$, $B_n$ with $n$ even, or
$D_n$ with $n$ odd, then $n_P = 2$.
\item"{(iii)}" In all other cases, $n_P =1$.
\endroster
Let $\bar G$ be the quotient of $G$ by the unique subgroup of
the center of $G$ of order $n_P$. Then the universal
$G$-bundle over $E\times H^1(E; \underline{U}(\eta))$ descends
to a universal
$\bar G$-bundle $\overline{\Xi}$ on $E\times \Cal M_E^0$.
\endstatement

Let us mention the analogous results for families of elliptic
curves over a base $B$. Let $\pi \: Z \to B$ be a flat family,
all of whose fibers are smooth elliptic curves or more
generally irreducible curves of arithmetic genus one (i\.e\.
smooth, nodal, or cuspidal curves). Let $\sigma$ be a section
of $\pi$ meeting each fiber in a smooth point. Associated to
$Z$ is the line bundle $L$ on $B$ defined by
$L^{-1} = R^1\pi _*\scrO_Z$, which can be identified with
$\scrO_Z(\sigma)|\sigma$ under the isomorphism $\sigma \to B$
induced by $\pi$. We want to describe the parabolic
construction along the family
$Z$. To do so, recall that we have the weights $g_i$ of (2.8),
which we assume ordered so that
$g_i \leq g_{i+1}$. Recall also that we have defined the
Casimir weights
$d_1, \dots, d_r$ of a root system $R$ in Section 3.  We order
the $d_i$ by increasing size, except in the case of $D_n$,
where we order the $d_i$ by: $2, 4,n,6,  8, \dots, 2n-2$.

Our result in families can then be somewhat loosely stated as
follows:

\theorem{4.5} Suppose that $G\neq E_8$. The parabolic
construction then globalizes over $Z$ to give a bundle of
nonabelian cohomology groups over
$B$. This bundle is a bundle of affine spaces with a
$\Cee^*$-action which is isomorphic to the vector bundle
$$\scrO_B \oplus L^{-d_1} \oplus \cdots \oplus L^{-d_r}.$$ Via
this isomorphism $\Cee ^*$ acts diagonally on the line bundles
in the direct sum, by the weight
$g_i$ on the factor $ L^{-d_i}$ \rom(and with weight $g_0 =
1$ on the factor
$\scrO_B$\rom). The associated bundle of weighted projective
spaces is then a universal relative moduli space for
$G$-bundles which are regular and semistable on every fiber.
\endstatement

A result closely related to Theorem 4.5 was established by
Wirthm\"uller \cite{33}, who also noted the exceptional
status of
$E_8$. We note that, from our point of view, in the case of
$E_8$ there is a family of weighted projective spaces over the
open subset $B'$ of $B$ over which the fibers of $\pi$ are
either smooth or nodal. However, this family is not the
quotient of a vector bundle minus its zero section by $\Cee^*$
acting diagonally. Furthermore, the construction degenerates
in an essential way at the cuspidal curves. A similar
phenomenon appears if we try to classify generalized del Pezzo
surfaces of degree one with an appropriate hyperplane section.

\section{5. Automorphism sheaves and spectral covers.} 

In this section, we fix  $G$ and denote by $\Cal M_E(G)=\Cal
M_E$ the moduli space of regular semistable $G$-bundles over
$E$. Likewise, given an elliptic fibration with a section $\pi
\: Z \to B$ whose fibers are smooth elliptic curves or nodal or
cuspidal cubics (except in the case $G=E_8$ where we will not
allow cuspidal fibers), we have a relative
moduli space $\Cal M_{Z/B}=
\Cal M_{Z/B}(G)$. Thus in all cases $\Cal M_{Z/B}$ is a
bundle of weighted projective spaces. 

Because of fixed points for the
$\Cee ^*$ action, the universal $G$-bundle over $E\times
H^1(E; 
\underline{U}(\eta))$ does not descend to a universal
$G$-bundle over $E\times \Cal M_E$, even locally, near the
singular points of
$\Cal M_E$, and a similar statement holds in families.
However, let
$\Cal M_E^0$ denote the smooth locus of $\Cal M_E$, and
similarly for $\Cal M_{Z/B}^0$. Then locally in either the
classical or
\'etale topology there exists a universal bundle $\Xi$ over
$E\times \Cal M_E^0$, and similarly for $Z\times _B\Cal
M_{Z/B}^0$. As we have seen in Theorem 4.4, there also exists
a $\bar G$-bundle $\overline{\Xi}$ over $E\times \Cal M_E^0$,
where $\bar G$ is a quotient of $G$ by a subgroup of the
center of order at most two. In particular, a universal
adjoint bundle always exists. In this section, we describe the
issues of the existence and uniqueness of a global universal
bundle over $E\times \Cal M_E^0$ or  $Z\times _B\Cal
M_{Z/B}^0$.

There are other questions closely related to these. Given a
family $\pi \: Z \to B$ as above, suppose that
$\Xi$ is a $G$-bundle over $Z$ such that $\Xi |\pi ^{-1}(b)$
is a semistable bundle for all $b$ for which $\pi ^{-1}(b)$ is
smooth. Then $\Xi$ defines a section of $\Cal M_{Z/B}$ over the
open subset of $B$ consisting of such $b$. At the singular
points of $\Cal M_{Z/B}$, the section is locally liftable to
the affine bundle of cohomology groups over $B$. Conversely,
a locally liftable section defines local
$G$-bundles over $\pi ^{-1}(U)$ for all sufficiently small
open sets $U$ of $B$ (in the classical or \'etale topology).
Note that the parabolic construction extends over the singular
fibers of $\pi$ (except for cuspidal fibers in case $G=E_8$),
dictating the correct definition of regular semistable
$G$-bundles for a singular fiber. When does a locally liftable
section of $\Cal M_{Z/B}$ actually determine a
$G$-bundle over
$Z$? More generally, how can we  describe the set (possibly
empty) of all bundles corresponding to a given section? For
simplicity, we shall assume that the section does not pass
through the singular points of
$\Cal M_{Z/B}$. Thus, if there existed a relative universal
bundle over $Z\times _B\Cal M_{Z/B}^0$, we could simply pull
this bundle back by the section to obtain a bundle over $B$.
While a relative universal bundle does not usually exist,
there are many cases where a section does indeed determine a
$G$-bundle. However, our answers are complete only in the
cases $G = SL(n), Sp(2n)$. 

Working for the moment with a single curve $E$, over an open
subset of $\Cal M_E^0$ where there exists a local universal
bundle $\Xi$, there is an associated group scheme
$\underline{\Aut}(\Xi)$. Because the associated automorphism
groups are abelian on $\Cal M_E^0$, as follows from (3.10),
these local group schemes piece together to give an abelian
group scheme over
$\Cal M_E^0$, whose associated sheaf of sections will be
denoted
$\Cal A$. In the usual way, the obstruction to finding a
global universal principal $G$-bundle over $E\times
\Cal M_E^0$ lies in  $H^2(\Cal M_E^0; \Cal A)$, and if this
obstruction is zero, then the set of all such principal
bundles is a principal homogeneous space over  $H^1(\Cal
M_E^0; \Cal A)$. More generally, given an elliptic fibration
$\pi\: Z \to B$ as above, we can fit together the automorphism
group schemes of local universal bundles to find an abelian
group scheme  over $\Cal M_{Z/B}^0$ whose fiber over every
point $b\in B$ is the group scheme constructed above. Let $\Cal
A_B$ denote the sheaf of sections of this group scheme. Given a
section $s$ of $\Cal M_{Z/B}^0
\to B$, we can pull back the the above group scheme to obtain
a group scheme over $B$, whose sheaf of sections we denote by
$\Cal A_B(s)$. Just as in the
case of a single smooth elliptic curve, the obstruction to 
finding a $G$-bundle over $Z$ corresponding to the section $s$
lies in
$H^2(B; \Cal A_B(s))$, and if this obstruction is zero, then
the set of all such  bundles is a principal homogeneous space
over  $H^1(B; \Cal A_B(s))$. 

Let us describe the sheaf $\Cal A$ in the case of $SL(n)$ and
a fixed elliptic curve $E$ in more detail. For each integer
$d$, $1\leq d
\leq n-1$, one can construct a universal extension $\Cal E_d$
over $E\times \Pee ^{n-1}$, viewing $\Pee^{n-1}$ as $\Ext
^1(W_{n-d}, W_d^*)$, which fits into an exact sequence
$$0 \to \pi _1^*W_d^* \otimes \pi _2^*\scrO_{\Pee ^{n-1}}(1)
\to \Cal E_d \to
\pi _1^*W_{n-d} \to 0.$$  Clearly, $\det \Cal E_d$ has trivial
restriction to each slice $E\times \{s\}$ but is not in fact
trivial. On the other hand, since the restriction of $\Cal
E_d$ to every fiber is regular and semistable,
$\pi _2{}_*Hom (\Cal E_d, \Cal E_d)$ is a sheaf of locally
free commutative
$\Cee$-algebras over $\Pee^{n-1}$ of rank $n$, and thus
corresponds to a finite morphism $\nu \: T \to \Pee ^{n-1}$ of
degree $n$, which we shall call the {\sl spectral cover\/} of
$\Pee ^{n-1}$. It is straightforward to identify the base
$\Pee ^{n-1}$ with the complete linear system $|np_0|$ and the
cover $T$ with the incidence correspondence in $\Pee
^{n-1}\times E$, in other words
$$T = \left\{\, (\sum _{i=1}^ne_i, e): \sum _{i=1}^ne_i \in
|np_0|, e = e_i {\text{ for some $i$}}\,\right\}.$$ Thus $T$
is smooth, and it has the structure of a $\Pee ^{n-2}$-bundle
over $E$ such that the $\Pee ^{n-2}$ fibers are mapped to
hyperplanes in $\Pee^{n-1}$ under $\nu$. Another way to
describe
$T$ is as follows: let
$\Lambda\cong \Zee ^{n-1}$ as the sublattice of $\Zee ^n$ of
vectors whose sum is zero, acted on by the Weyl group $\frak
S_n$, so that $\Pee ^{n-1} = |np_0| = (E\otimes \Lambda)/\frak
S_n$. Let
$W_0 =\frak S_{n-1}\subset \frak S_n$ be the stabilizer of the
vector
$e_n \in \Zee ^n$. Then $T = (E\otimes \Lambda)/W_0$.

A standard argument shows that, if
$\Cal V$ is a vector bundle  over $E\times \Pee^{n-1}$ whose
restriction to every slice is isomorphic to the corresponding
restriction of $\Cal E_d$, then
$\pi _2{}_*Hom (\Cal V, \Cal E_d)$ is locally free of rank one
over $\pi _2{}_*Hom (\Cal E_d,
\Cal E_d) = \nu _*\scrO_T$, and thus corresponds to a line
bundle on $T$, and conversely every line bundle on $T$ defines
a vector bundle $\Cal V$ with the above property. It is
helpful to compare this situation with the one usually
encountered in algebraic geometry, where we try to make a
moduli space of simple vector bundles and then the only choice
is to twist by the pullback of a line bundle from the moduli
space factor. 

From this it follows that, in the case of $SL(n)$, the
automorphism sheaf $\Cal A$ is given by the kernel of the norm
homomorphism $\nu_*\scrO_T^* \to \scrO_{\Pee ^{n-1}}^*$. Hence
there is an exact sequence
$$0 \to H^1(\Pee^{n-1}; \Cal A) \to \Pic T \to \Pic \Pee
^{n-1} \to H^2(\Pee^{n-1}; \Cal A) \to H^2(\scrO_T^*)\to 0.$$
Thus, 
$H^1(\Pee^{n-1}; \Cal A) \cong \Zee\times E$ for $n>2$ and
$H^1(\Pee^1; \Cal A) \cong E$, and
$H^2(\Pee^{n-1}; \Cal A)
\cong H^3(T; \Zee)$. There is also an analogue of the above
exact sequence where we take \'etale cohomology. In this case,
$H^1(\Pee^{n-1}; \Cal A)$ is unchanged and $H^2(\Pee^{n-1};
\Cal A) \cong H_{\text{\'et}}^2(T; \Bbb G_m)$, which is a
torsion group.  The obstruction to gluing together local
families (in either the classical or \'etale topology) of
$SL(n)$-bundles to make a global
$SL(n)$-bundle over
$E\times \Pee^{n-1}$ lives in $H^2(\Pee^{n-1}; \Cal A)$, and in
case the obstruction is zero the set of all such bundles in
then a principal homogeneous space over
$H^1(\Pee^{n-1}; \Cal A)$. In our case, a direct construction
using the pushforward of appropriate line bundles on $E\times
T$ shows that the obstruction in $H^2(\Pee^{n-1};
\Cal A)$ vanishes. Thus the  family of universal 
$SL(n)$-bundles
$\Cal V$ over $E\times \Pee^{n-1}$ is parametrized by
$\Zee\times E$ for $n>2$ and by $E$ if $n=2$. In case we
consider the corresponding situation in families $Z\to B$,
then there  exist mod
$2$ obstructions to finding a principal $SL(n)$-bundle over the
entire family, and these obstructions are not in general zero.
On the other hand, there always exists a universal
$GL(n)$-bundle
$V$ such that
$V|\pi ^{-1}(b)$ has trivial determinant for all $b$, so that
$\det V$ is pulled back from $B$. See
\cite{20} for more detail in the case of vector bundles.

Similar explicit constructions can be carried out for the
symplectic group.  Let $\Lambda = \Zee ^n$, and let the Weyl
group $W = \frak S_n \ltimes (\Zee/2\Zee)^n$ act on $\Lambda$,
where the symmetric group acts by permuting the basis elements
and $(\Zee/2\Zee)^n$ acts by sign changes. Then $(E\otimes
\Lambda)/W = \Sym ^n\Pee ^1 = \Pee ^n$. Let $W_0 = \frak
S_{n-1}
\ltimes (\Zee/2\Zee)^{n-1}$ be the subgroup of $W$ fixing the
last basis vector, and set $T^{\text{sp}} = (E\otimes
\Lambda)/W_0 = \Pee ^{n-1} \times T$. The group $W_0$ is a
subgroup of index two in the larger group $W_1 = \frak S_{n-1}
\ltimes (\Zee/2\Zee)^n$, and there is an induced involution
$\iota$ on $T^{\text{sp}}$, with quotient $T^{\text{sp}}/\iota
= S =  \Pee ^{n-1} \times \Pee ^1$. We then have:

\proposition{5.1} For the symplectic group $Sp(2n)$, the
automorphism sheaf $\Cal A^{\text{sp}}$ over $\Cal M_E(Sp(2n))
\cong \Pee^n$  is given by
$$\Cal A^{\text{sp}} = \{\, f\in \nu_*\scrO_{T^{\text{sp}}}^*:
\iota ^*f= f^{-1}\,\}.$$
\endproclaim

Using (5.1) one can show that there is a universal bundle over
$E\times
\Cal M_E(Sp(2n))$ as well, and that the set of all universal
bundles is parametrized by $E$. Thus we have constructed
universal bundles over $E\times \Cal M_E$ in the two cases
where the moduli space is smooth. It then follows from Theorem
4.4 that a universal bundle exists over
$\Cal M_E^0(G)$ in all cases, with the possible exception of
$G= Spin (4n+1)$ and $G = Spin (4n+2)$.

We return to the  case of a general $G$ and analyze the
structure of the sheaf $\Cal A$ over $\Cal M_E^0$. Since
$\Cal A$ is the sheaf of sections of an abelian algebraic
group scheme, there is the exponential map $\exp$ from the
corresponding sheaf of Lie algebras
$\operatorname{Lie}\Cal A$ to $\Cal A$. The kernel of $\exp$
is a constructible sheaf, which we denote by
$\underline{\Lambda}$, and the image of $\exp$ is the sheaf
$\Cal A^0$ which, locally, consists of all sections of
$\underline{\Aut}(\Xi)$ passing through the identity component
of every fiber. First we note that $\Cal A = \Cal A^0$ on the
Zariski open subset $U$ of $\Cal M^0$ consisting of split
bundles, where the fiber over $x\in U$ of the group scheme
corresponding to $\Cal A$ is $(\Cee ^*)^r$ and is connected. If
the root system for
$G$ is simply laced, we can say more:

\proposition{5.2} Suppose that $G$ is simply laced.  If $G\neq
SL(2)$, then the set
$$\{\, \xi \in \Cal M^0: \Cal A_\xi \neq \Cal A^0_\xi\,\}$$
has codimension at least two in $\Cal M^0$. 
\endproclaim

As a consequence, in the relative setting, for $G$ simply
laced, if
$\dim B = 1$ and $G \neq SL(2)$, then for a generic section
$s$ of
$\Cal M_{Z/B}^0$, we can always assume that $\Cal A_B(s) = \Cal
A_B^0(s)$. The above
proposition does not hold if $G$ is not simply laced; for
example, it fails for
$Sp(2n)$.

Next we turn to $\Cal A^0 =\operatorname{Lie}\Cal A
/\underline{\Lambda}$. Note that, in case there is a universal
bundle $\Xi$ over $E\times \Cal M^0$, then
$\operatorname{Lie}\Cal A = R^0p_2{}_*(\ad \Xi)$ is dual to
$R^1p_2{}_*(\ad \Xi)$, which is the tangent bundle to $\Cal
M^0$.  Thus
$\operatorname{Lie}\Cal A =
\Omega ^1_{\Cal M^0}$ is the cotangent bundle. In fact, this
statement always holds, since a universal bundle exists
locally and the automorphism sheaf is abelian.  Another way to
describe the cotangent bundle is as follows: let $\frak h$ be
the Lie algebra of a Cartan subgroup of $G$. Then the Weyl
group acts on $E\otimes
\Lambda$ and on the trivial vector bundle $\scrO_{E\otimes
\Lambda}\otimes _\Cee\frak h$, and the sheaf of $W$-invariant
sections is a coherent sheaf over $(E\otimes \Lambda)/W= \Cal
M$ whose restriction to $\Cal M^0$ is locally free, and in
fact is
$\Omega ^1_{\Cal M^0}$. The constructible sheaf
$\underline{\Lambda}$ can be described as follows. Let $U$ be
the open subset of $\Cal M^0$ over which the map $E\otimes
\Lambda \to
\Cal M$ is unramified, and let $i\: U \to \Cal M^0$ be the
inclusion. Then the action of $W$ on $\Lambda$ gives a locally
constant sheaf $\underline{\Lambda} _0$ on $U$, and
$\underline{\Lambda} = i_*\underline{\Lambda}_0$. The map
$\Lambda \to \frak h$ induces an inclusion
$\underline{\Lambda} \to \left(\scrO_{E\otimes
\Lambda}\otimes _\Cee\frak h\right)^W$, and this is the same
as the inclusion $\underline{\Lambda} \to
\operatorname{Lie}\Cal A$.

This picture is related to the general theory of spectral
covers of \cite{24} and \cite{10} (as has also been
noted by Donagi in
\cite{11}). Suppose that
$\varpi$ is an element of $\frak h$ such that $W\cdot\varpi$
spans
$\frak h$ over
$\Cee$. In the typical application, $\varpi$ is (the dual of) a
minuscule weight, if such exist. Let $W_0$ be the stabilizer of
$\varpi$. If we set $T = (E\otimes \Lambda)/W_0$, then there
is a surjection $\nu \: T \to \Cal M$. By pure algebra,
$$\nu _*\scrO_T = \left(\scrO_{E\otimes \Lambda}\otimes_\Cee
\Cee[W/W_0]\right)^W.$$ On the other hand, there is a
surjection $\Cee[W/W_0] \to \frak h$ whose kernel consists of
the relations in the orbit $W\cdot
\varpi$. Correspondingly, there is a surjection
$$\left(\scrO_{E\otimes \Lambda}\otimes_\Cee
\Cee[W/W_0]\right)^W \to \left(\scrO_{E\otimes
\Lambda}\otimes_\Cee
\frak h\right)^W.$$ In particular $H^1(\Cal
M;\operatorname{Lie}
\Cal A)$ is a quotient of
$H^1(\Cal M; \nu_*\scrO_T)$.

Now suppose that we are in the relative case of an elliptic
fibration $\pi \: Z\to B$. There is then a relative universal
moduli space $\Cal M_{Z/B}$ (with the usual care in the case of
$E_8$). The covers $T \to \Cal M$ defined over every smooth
fiber extend to a finite morphism $\Cal T_{Z/B} \to \Cal
M_{Z/B}$. A section $s$ of the map $\Cal M^0_{Z/B} \to B$
defines a finite cover
$C_s$ of $B$, which we will call the {\sl spectral cover\/} in
this case. Of course, $C_s$ need not be smooth or even
reduced. In case
$\dim B = 1$ and $s$ is generic, the above discussion
identifies the connected components of
$H^1(B;\Cal A_B(s))$ with an abelian variety which is a
quotient of the Jacobian $J(C_s)$, and which is called the {\sl
Prym-Tyurin variety\/} of the spectral cover.

A straightforward dimension argument shows:

\proposition{5.3} Suppose that $\dim B = 1$ and that $\Cal
A_B(s)_b = \Cal A_B^0(s)_b$ for at least one point $b\in B$.
Then $H^2(B; \Cal A_B(s)) = 0$. In other words, there exists a
universal $G$-bundle over $B$ corresponding to the section $s$.
\endstatement

If however $\Cal A_B(s)_b \neq \Cal A_B^0(s)_b$ for all $b\in
B$, then it is possible for there not to exist a universal
$G$-bundle over $B$ corresponding to $s$, even when $G=
SL(2)$. For $\dim B$ arbitrary, the possible obstructions in
the case of $SL(n)$ are analyzed in detail in \cite{20}. 

Let us work out the twisting group $H^1(B;\Cal A_B(s))$
explicitly in the simplest cases $G= SL(n), Sp(2n)$, with
$\dim B$ arbitrary:

\proposition{5.4} Suppose that $G = SL(n)$. Let $C_s \to B$ be
the spectral cover defined above. Then
$$H^1(B;\Cal A_B(s)) = \Ker \{\, \operatorname{Norm}\: \Pic
(C_s) \to
\Pic B\,\}.$$ If $G = Sp(2n)$, let $C_s$ be the corresponding
degree $2n$ cover of $B$, let $\iota\: C_s\to C_s$ be the
induced involution, and let
$f\: C_s\to D_s$ be the degree two quotient of $C_s$ by
$\iota$. Then 
$$H^1(B;\Cal A_B(s)) = \Ker \{\, \operatorname{Norm}\: \Pic
(C_s) \to
\Pic D_s\,\}.$$
\endstatement

Thus, in case $B$ is a curve, $H^1(B;\Cal A_B(s))$ is the
generalized Prym variety of the cover $C_s\to D_s$. Similar
results hold for the remaining classical groups $Spin (2n)$ and
$Spin (2n+1)$. 

On the other hand, suppose that $G= E_6, E_7, E_8$, that $\dim
B = 1$, and that the section $s$ is generic. In this case,
there is an associated fibration of del Pezzo surfaces $p\: Y
\to B$, where $Y$ is a smooth threefold. Moreover $Z$ is
included as a smooth divisor on $Y$ so that $p|Z = \pi \: Z
\to B$. Let
$J^3(Y)$ denote the intermediate Jacobian of $Y$. There is an
induced morphism
$J^3(Y) \to J(B)$, where $J(B)$ is the ordinary Jacobian of
$B$ coming from the homomorphism $H^*(Y) \to H^*(Z) \to
H^{*-2}(B)$.  Denote the kernel of  the morphism $J^3(Y) \to
J(B)$ by
$J^3(Y/B)$. Finally set
$$H^{2,2}_0(Y; \Zee) = \left.\Ker \{\, H^4(Y; \Zee) \to H^2(B;
\Zee)\,\}\right/\Zee\cdot [Y_t],$$ where $Y_t$ is a general
fiber of $p$. In general $H^{2,2}_0(Y;
\Zee)$ is a finite group.  We then obtain the following
theorem, first proved by Kanev
\cite{24} in the case
$B =\Pee^1$ via the Abel-Jacobi homomorphism:

\proposition{5.5} In the above situation, there is an exact
sequence
$$0\to J^3(Y/B) \to H^1(\Cal A_B(s)) \to H^{2,2}_0(Y; \Zee)
\to 0.$$
\endstatement

We note that we can interpret $H^1(\Cal A_B(s))$ as a relative
Deligne cohomology group.

\enddocument
\end